%% file: main.tex
\documentclass[letterpaper, 10 pt, conference]{ieeeconf}
\IEEEoverridecommandlockouts

\usepackage{amsmath,amsfonts,amssymb}

\usepackage{graphics} 
\usepackage{epsfig,color} 
\usepackage{times} 
\usepackage{amssymb}
\usepackage{wrapfig}
\usepackage{amsfonts}
\usepackage{amsthm}
\usepackage{graphicx}
\graphicspath{ {./images/} }
\usepackage{graphicx}
\usepackage{subfig}
\usepackage{amsmath}
\usepackage[ruled,vlined]{algorithm2e}
\usepackage{footnote}
\newtheorem{assumption}{Assumption}
\newtheorem{remark}{Remark}
\newtheorem{theorem}{Theorem}

\newtheorem{lemma}[theorem]{Lemma}
\newtheorem{definition}{Definition}

\title{\LARGE \bf
Data-Driven Near-Optimal Control of Nonlinear Systems Over Finite Horizon
}
\author{Vasanth Reddy, Hoda Eldardiry and Almuatazbellah Boker
\thanks{Vasanth Reddy and Hoda Eldardiry are with the Department of Computer Science,
        Virginia Tech, Blacksburg, VA, USA.
        {\tt\small emails: vasanth2608@vt.edu and hdardiry@vt.edu}}%
\thanks{Almuatazbellah Boker is with Bradley Department of Electrical and Computer Engineering, Virginia Tech, Blacksburg, VA, USA.
        {\tt\small email: boker@vt.edu}}%
}

\begin{document}

\maketitle

\begin{abstract}
We examine the problem of two-point boundary optimal control of nonlinear systems over finite-horizon time periods with unknown model dynamics by employing reinforcement learning. We use techniques from singular perturbation theory to decompose the control problem over the finite horizon into two sub-problems, each solved over an infinite horizon. In the process, we avoid the need to solve the time-varying Hamilton-Jacobi-Bellman equation. Using a policy iteration method, which is made feasible as a result of this decomposition, it is now possible to learn the controller gains of both sub-problems. The overall control is then formed by piecing together the solutions to the two sub-problems. We show that the performance of the proposed closed-loop system approaches that of the model-based optimal performance as the time horizon gets long. Finally, we provide three simulation scenarios to support the paper's claims. 
\end{abstract}


\IEEEpeerreviewmaketitle

\input{Sections/1.Introduction}
\input{Sections/2.Problem_formulation}
\input{Sections/3.Method}
\input{Sections/4.Learning}
\input{Sections/5.Simulation}
\input{Sections/6.Conclusion}
\bibliographystyle{refstyle}
\bibliography{main}


\end{document}

%% file: Sections/1.Introduction.tex
\section{Introduction}
Finding the optimal controller for a nonlinear system with a finite time horizon is significantly more challenging than finding the optimal controller for nonlinear systems with an infinite time horizon. This is due to the fact that finding the optimal controller with a finite time horizon requires the solution of a time-varying HJB equation. In addition, the dynamics may be unknown or susceptible to modeling uncertainties. To address these challenges, researchers have developed a wide variety of reinforcement learning-based techniques. Previous works have employed adaptive dynamic programming and the actor-critic learning strategy to train the controller for nonlinear systems, as documented in research by~\cite{jiang2015global, lv2016online} and~\cite{vamvoudakis2017q, vamvoudakis2010online, luo2016model, lin2017online}, respectively. Nonetheless, each of the aforementioned methods focused on problems involving optimal control over infinite horizons. On the other hand, very few learning algorithms help in learning the controller for nonlinear systems across finite horizons. In this respect, a few model-based~\cite {kim2018deep} and model-free~\cite{zhao2015finite, zhao2020finite, chen2022two} approaches learn the weights to the time-varying basis function vectors to produce a solution. What makes this problem particularly challenging is the design of time-varying basis functions. We, on the other hand, avoid this challenge by transforming the HJB equation into a time-invariant one by employing a singular perturbation method. 

In this study, we solve a two-point boundary optimal control problem for an uncertain nonlinear system over a finite time horizon by learning the controller. Key to the proposed solution is the recognition that the two-point boundary optimal control problem exhibits a two-time scale phenomenon, even when the original system is not singularly perturbed, as demonstrated in the seminal work~\cite{wilde1972dichotomy}. In this case, the dynamics of the closed-loop system evolve at a faster rate in comparison to the time horizon. This phenomenon becomes more apparent as the horizon time gets long. This two-time-scale phenomenon allows the optimal control problem over finite horizons to be approximated by two sub-problems over infinite horizons. More specifically, the time decomposition of the problem allows for casting the problem into stabilizing the system forward in time and, respectively, backward in time then attaching the two solutions together. This results in a near-optimal performance, which gets closer to the optimal one as the time horizon gets long. Furthermore, we utilize a policy iteration technique via continuous time Q-learning, recently developed in~\cite{chen2019adaptive}, to estimate the control gains of the two boundary sub-problems when considering the uncertainty of the model dynamics. The combination of the time decomposition of the problem and utilization of data for controller learning leads to the main contribution of the paper. That is, we solve the two-point boundary optimal control problem for a large class of uncertain nonlinear systems using policy iteration methods without the need to solve time-varying or partial deferential equations. This results in a relatively simple and efficient learning algorithm.   

The following outline constitutes the framework of this study.  Section~\ref{problem formulation} includes a schematic that explains the system setup as well as the formulation of the problem. The section~\ref{method} details the original problem's two-time scale reduction. The policy iteration learning procedure that is used to estimate the system's control gains is described in Section~\ref{learning}. Section~\ref{simulation} has a simulation example and Section~\ref{conclusion} has final remarks.

%% file: Sections/2.Problem_formulation.tex
\section{Problem Formulation} \label{problem formulation}
Consider a continuous-time nonlinear affine system of the form
\begin{equation} \label{system}
    \dot{x} = f(x) + g(x)u(x), \quad x(0) = x_0,\; x(T) = x_T
\end{equation}
where $x \in \mathbb{R}^n$ and $u(x) \in \mathbb{R}^m$ are the system states, control input, respectively. 
\begin{assumption}
The nonlinear functions $f(.)$ and $g(.)$ are smooth, and globally Lipschitz. In addition, $f(.)$ is unknown. 
\end{assumption}
\begin{assumption}
     The system~\eqref{system} is controllable from any initial state to the origin and from the origin to any prescribed state. 
\end{assumption}
It is the aim to design the control law $u(x)$ that will drive the system states, $x(t)$ from the beginning state $x(0) = x_0$ to the final terminal state $x(T) = x_T$ during a time period $T$ and, at the same time, minimize the value of the objective function.
\begin{equation} \label{objectivefunction}
    J=  \int_{0}^T \mathcal{R}(x(\bar{t}), u(\bar{t})) \quad d \bar{t}
\end{equation}
where $\mathcal{R}(x(t), u(t)) = \mathcal{S}(x(t)) + u^\top(t) R u(t)$, $\mathcal{S}(t) \succ 0$ and $R \succ 0\;\;\forall\; t \in [0, T].$ 
We define the value function $V(x(t))$ as:
\begin{equation} \label{valuefunction}
    V^*(x(0), 0, T) = \min_{u^*(t)} J(x, u^*(x), 0, T),
\end{equation}
where $V^{*}$ is Bellman's optimal value function and $u^*(x)$ is the optimal control. This optimal control problem can be solved by solving the Hamiltonian-Jacobi-Bellman equation  \cite{athans2013optimal} 
\begin{equation} \label{HJB}
    -\frac{\partial V^{*}}{\partial t}=\min _{u^{*}(.)}\left[\mathcal{R}(x, u)+\frac{\partial V^{*}}{\partial x} (f(x) + g(x)u^{*}(x))\right].
\end{equation}
 The HJB equation \eqref{HJB} is a partial differential equation that describes the optimal control policy. However, it is often too complex to solve analytically. In the following section, we will follow a singular perturbation approach to simplify~\eqref{HJB} in the limit when $T$ is relatively long. We will then use this result to 
 provide a data-driven method to learn the controller. 


%% file: Sections/3.Method.tex
\section{Control Design using Singular-Perturbation analysis} \label{method}
For control problems that are defined over a finite time period, it is known that when $T$ gets large, the trajectories of the system start displaying two time-scale features \cite{wilde1972dichotomy}. In this case, the majority of transient activity takes place at the boundary points where as the dynamics remain in a steady state in between. 
Because of this phenomenon, it is possible to concentrate on the way the system behaves around the points $t=0$ and $t=T$. This will lead to a simplification of the control problem as will be shown next.


We start by setting up the singular perturbation model of the system. For this purpose, we normalize the time period $0$ to $T$ to the interval $[0, 1]$ and introduce the scaled time parameter $\varepsilon$ as 
\begin{equation}
    \tau = \frac{t}{T} ,\quad \varepsilon =1/T.   \label{tau}
\end{equation}
In view of \eqref{tau}, we rewrite~\eqref{system},~\eqref{objectivefunction} and~\eqref{HJB} to obtain
\begin{align}
    \varepsilon \dot{x}(\tau) &= f(x(\tau)) + g(x(\tau))u(\tau), \quad x(0) = x_0, x(1) = x_T, \label{scaledsystem}\\
    -\varepsilon \frac{\partial V^{*}}{\partial \tau}&=\min _{u^{*}(\tau)}\left[\mathcal{R}(x, u)+\frac{\partial V^{*}}{\partial x} (f(x) + g(x)u^{*}(\tau))\right]. \label{scaledhjb}
\end{align} 
Letting $\varepsilon \to 0$ in \eqref{scaledhjb} leads to
\begin{equation}
    0=\min _{u^{*}(\tau)}\left[\mathcal{R}(x, u^*)+\frac{\partial V^{*}}{\partial x} (f(x) + g(x)u^{*}(\tau))\right]. \label{hjbsp}
\end{equation}
Following a singular perturbation approach, it can be shown that \eqref{hjbsp} has two boundary layer solutions $u_{+}(x)$ and $u_{-}(x)$ \cite{kokotovic1999singular}.  
As in the case of linear systems~\cite{reddy2022singular}, the \textit{forward controller} $u_{+}(x)$ will be designed to stabilize the system~\eqref{scaledsystem} in real time and the \textit{backward controller} $u_{-}(x)$ will be designed to stabilize the system in reverse  time. Next, we will explicitly describe these two regulators.

\subsection{Forward Regulator}

In then limit as $\varepsilon \to 0$, the value function for the system~\eqref{scaledsystem} is given as:
\begin{equation} \label{valuefunction+}
    V_{+}(x_+(0)) = \min_{u_+(.)} \int_{0}^\infty \mathcal{R}(x_+(\bar{\tau}), u_+(\bar{x})) \quad d \bar{\tau}.
\end{equation}
The associated value HJB optimality condition for the value function~\eqref{valuefunction+} is given as: 
\begin{equation} \label{HJB+}
    \min_{u_+(.)}\left[\mathcal{R}(x_+, u_+)+\frac{dV_{+}^{*}}{dx_+} (f(x_+) + g(x_+)u_+^{*}(x_+)))\right] = 0 
\end{equation}
Assuming that the above minimum exists and is unique, then the optimal control function for the given system~\eqref{scaledsystem} is provided as: 
\begin{equation} \label{control+}
    u^{*}_{+}(x) = -\frac{1}{2}R^{-1}g^\top(x_+)\frac{dV_{+}^{*}}{dx_+}.
\end{equation}
\subsection{Backward Regulator}
Now consider the reverse time $s = -\tau $, then 
the associated value function for the system $\frac{dx}{ds} = - f(x)  - g(x)u$ is given as:
\begin{equation} \label{valuefunction-}
    V_{-}(x_-(T)) = - \int_{-\infty}^T \mathcal{R}(x_-(\bar{s}), u_-(\bar{s})) \quad d \bar{s}
\end{equation}

The associated value HJB optimality condition for the value function~\eqref{valuefunction-} is given as
\begin{equation} \label{HJB-}
    \min_{u_-(.)}\left[\mathcal{R}(x_-, u_-)-\frac{dV_{-}^{*}}{dx_-} (f(x_-) + g(x_-)u_-^{*}(x_-)))\right] = 0 
\end{equation}
Assuming that the above minimum exists and is unique, then the optimal control policy is provided as: 
\begin{equation} \label{control-}
    u^{*}_{-}(x) = \frac{1}{2}R^{-1}g^\top(x_-)\frac{dV_{-}^{*}}{dx_-}.
\end{equation}\
\subsection{Near Optimal Performance}
It is shown in Property 4.2 in~\cite{anderson1987optimal} and Theorem 2.1 Chapter 6 in~\cite{kokotovic1999singular} that the combination of the two infinite time horizon problems will approximate that of the original value function~\eqref{valuefunction} of system~\eqref{system} for sufficiently small $\varepsilon$. The following theorem summarizes this result. 

\begin{theorem} \label{th2}
Let Assumptions 1-2 hold. Suppose further that $\|x(0)\| \leq 1$ and $\|x(T)\| \leq 1$. Then there exists $\varepsilon_1>0$ such that, for all $\varepsilon\in(0,\varepsilon_1]$,   
\begin{align}
    V^*(x(0), 0, T) &\leq V^*_+(x(0)) - V^*_-(x(T)) + k_1(\varepsilon), \label{v1th}\\ 
    V^*_+(x(0)) -& V^*_-(x(T)) \leq V^*(x(0), 0, T) + k_2(\varepsilon), \label{v2th}
\end{align}
where $k_1(\varepsilon)$ and $k_2(\varepsilon)$ are monotonic functions of $\varepsilon$ with $\lim_{\varepsilon \to 0} k_1(\varepsilon) = 0$ and $\lim_{\varepsilon \to 0} k_2(\varepsilon) = 0.$ 
\end{theorem}

\begin{remark}
Together,~\eqref{v1th} and~\eqref{v2th} show that $V^*[x(0), x(T)]$ and $V^*_{+}[x(0)]-V^*_{-}[x(T)]$ can be made as close as desired by choosing $T$ large enough. 
\end{remark}
\begin{remark}
    In theorem~\ref{th2}, it is assumed that, $\|x(0)\| \leq 1$ and $\|x(T)\| \leq 1$. Nonetheless, any compact sets can be assumed for $x(0), x(T)$. The only effect of changing the sets is that the size of the interval $T$ to achieve a given closeness of $V^*[x(0), x(T)]$ and $V^*_{+}[x(0)]-V^*_{-}[x(T)]$ will be affected.
\end{remark}
In the next section, we develop a reinforcement learning algorithm to solve for the optimal value functions function~\eqref{valuefunction+} and ~\eqref{valuefunction-}, and hence the controllers~\eqref{control+} and~\eqref{control-}, without the need for $f(x)$ to be known.

%% file: Sections/4.Learning.tex
\section{Main Results} \label{learning}
\subsection{Learning Based Design}
 In this section, we follow a policy iteration method inspired by~\cite{chen2019adaptive} to learn the value functions of the forward and backward regulators. 
\subsubsection{Learning of the Forward Regulator}
 The key to our learning approach is Policy Iteration, which is the process of successive iterations of policy evaluation and improvement to arrive at the best possible policy. We will discuss policy iteration for the forward regulator problem next.\par 
\textbf{\emph{Policy Iteration: }}
\begin{enumerate}
    \item For a given control policy, $u_{+}^{k}(x_+)$, solve for the value function $V_{+}^{k}(x_+)$ using
    \begin{align}
        V_+^{k-1}(x_+(t-T))=\int_{t-T}^t \mathcal{R} (x_+(\tau), &u(x_+)) \mathrm{d} \tau \nonumber \\
        &+V_+^k(x_+(t)) \label{pi1}
    \end{align} 
    \item Revise the control policy using the updated value function
    \begin{equation} \label{pi2}
        u_{+}^{k+1} = -\frac{1}{2}R^{-1}g^\top(x_+)\nabla V_{+}^{k}(x)
    \end{equation}
\end{enumerate}
\par 
\textbf{\emph{Value Function Approximation: }} \par 
We now design the adaptive critic for policy evaluation. In this paper, we consider a neural network to approximate the value function and is given as: 
\begin{equation} \label{adaptivecritic}
    V_{+} = W_+^\top \phi (x_+) + \epsilon(x_+),
\end{equation}
where $\phi(x_+): \mathbb{R}^n \in \mathbb{R}^N$ represents the activation function vector
with the number $N$ of neurons in the hidden layer, $W_+ \in  \mathbb{R}^N$ represents the weight vector and $\epsilon(x_+) \in \mathbb{R}$ represents the neural network
approximation error. The activation functions are chosen to form a completely independent basis set, meaning that they can approximate any function uniformly. This choice ensures that the neural network can approximate any function $V_+(x_+)$ uniformly within a compact set $\Omega$. \par 
To update the adaptive critic using the Bellman approach, substitute~\eqref{adaptivecritic} in~\eqref{pi1} to obtain
\begin{equation} \label{irldiff}
   - \epsilon_{B} = \int_{t-T}^t \mathcal{R} (x_+(\tau), u(x_+)) \mathrm{d} \tau + W_{+}^{\top} \Delta \phi(t)
\end{equation}
where $\epsilon_{B} = \epsilon(x_+(T)) - \epsilon(x_+(t-T))$ is the Bellman residual error, $\Delta \phi(t) = \phi(x_+(t))- \phi(x_+(t-T)$ and defining the term, $\rho_+(x_+, u) = \int_{t-T}^t \mathcal{R} (x_+(\tau), u(x_+)) \mathrm{d} \tau$. \par 
We introduce a set of auxiliary variables called $\xi \in \mathbb{R}^
{N\times N}$ and $\psi \in \mathbb{R}^
{N}$ by low-pass filtering the variables in~\eqref{irldiff}. This is done so that we can construct an adaptive law that can estimate the weights of the value function approximation with guaranteed convergence~\cite{chen2019adaptive}. As a result, we have the adaptive law
\begin{align}
    \dot{\xi}&=-\ell \xi+\Delta \phi(t) \Delta \phi(t)^{\top}, \xi(0)=0, \label{aux1}\\ 
    \dot{\psi}&=-\ell \psi+\Delta \phi(t) \rho_{+}(x_+, u), \psi(0)=0,
\end{align}
where $\ell > 0$ ensures the parameters $\xi$ and $\psi$ to be bounded and guarantees stability~\cite{na2015robust}.
To update the weights $W_{+}$ of the critic network, we employ the sliding mode technique as given in~\cite{zhao2020finite}
\begin{equation} \label{weightupdate}
    \hat{W}_{+} = -\Gamma \xi \frac{G}{\|G\|}
\end{equation}
where $G = \xi \hat{W}_{+} + \psi$ and $\Gamma \succ 0$ is a learning gain parameter. 
The estimated adaptive critic function is given by:
\begin{equation} \label{criticest}
    \hat{V}_{+} = \hat{W}_+^\top \phi (x_+).
\end{equation}
 \par 
\textbf{\emph{Policy Evaluation: }} \par 
Now, we are going to construct an actor for the improvement of policy. If the weight $\hat{W}_+$ converges to the true unknown weight $W_+$, which solves the Bellman equation, one may discover the optimal control directly by inspecting~\eqref{pi2} and utilizing the adaptive critic~\eqref{criticest}. This is possible only if the Bellman equation is satisfied. The control policy will now be
\begin{equation}
    u_+=-\frac{1}{2} R^{-1} g(x_+)^{\top} \nabla \phi^{\top} \hat{W}_+.
\end{equation}

\subsubsection{Backward Regulator}
We follow the same steps as in the forward regulator problem except for the change in the integral term, $\rho_-(x_-, u) = -\int_{t-T}^t \mathcal{R} (x_-(\tau), u(x_-)) \mathrm{d} \tau$. The control policy for backward regulator problem will be $u_-=\frac{1}{2} R^{-1} g(x_-)^{\top} \nabla \phi^{\top} \hat{W}_-$. \par 

The pseudocode for the learning algorithm for both regulators is given in Algorithm~\ref{algo}. 
\begin{algorithm}[ht]
\SetAlgoLined
\KwResult{Control policies for forward and backward regulator problems}
Initialization: $W_+$, $W_-$, $\epsilon > 0,  k = 0$, exploration noise, $e$\\
\While{$\Delta W_+ \geq \epsilon$}{
Get the state, $x$ from the system~\eqref{system} using the control policy~\eqref{control+}, $u_+^k = u_+^k + e$ \\ 
Evaluate the value function $V_+^k = (\hat{W}_+^k)^\top \phi (x)$ \\ 
Policy evaluation: $u_+^{k+1} = -\frac{1}{2}R^{-1}g^\top(x)\nabla V_{+}^{k}(x)$ \\ 
Weight update as given in~\eqref{weightupdate}
Compute the change in weights: $\Delta W_+ = W_+^k  - W_+^{k+1}$;
$k = k + 1$;
}
$u_+ = -\frac{1}{2}R^{-1}g^\top(x)\nabla \phi^{\top} W_+^{k+1}$ \\
\While{$\Delta W_- \geq \epsilon$}{
Get the state, $x$ from the system~\eqref{system} using the control policy~\eqref{control-}, $u_-^k = u_-^k + e$ \\ 
Evaluate the value function $V_-^k = (\hat{W}_-^k)^\top \phi (x)$ \\ 
Policy evaluation: $u_-^{k+1} = \frac{1}{2}R^{-1}g^\top(x)\nabla V_{-}^{k}(x)$ \\ 
Weight update as given in~\eqref{weightupdate}
Compute the change in weights: $\Delta W_- = W_-^k  - W_-^{k+1}$;
$k = k + 1$;
}
$u_- = \frac{1}{2}R^{-1}g^\top(x)\nabla \phi^{\top} W_-^{k+1}$
\caption{Policy Iteration on two value boundary problems}
\label{algo}
\end{algorithm}

\subsection{Convergence Guarantees}
To guarantee parameter convergence of the adaptive function, we require a Persistent Excitation (PE) condition~\cite{chen2019adaptive}, which we define next.

\begin{definition} {Persistent Excitation (PE)~\cite{na2015robust}}
The signal
$\Delta \phi(x_+)$ is said to be persistently excited over the time interval
$[\tau-1, \tau]$ if there exists a strictly positive constant $\sigma_1 > 0$ such that
\begin{equation}
    \int_{\tau-1}^\tau \Delta \phi(\bar{\tau}) \Delta \phi(\bar{\tau})^{\top} \mathrm{d} \tau \geq \sigma_1 I, \quad \forall \tau>0.
\end{equation}
\end{definition}
We further need the following Lemma.
\begin{lemma}~\cite{na2015robust} \label{lemma1}
If the signal $\Delta \phi(x_+)$ is persistently excited for
all $\tau >$ 0, the auxiliary variable $\xi$ defined in~\eqref{aux1} is positive
definite, i.e. $\xi \succ 0$ and the minimum eigenvalue $\mathcal{\lambda}_{min}(\xi) >
\sigma_1 > 0, \forall \tau > 0$ for some positive constant $\sigma_1$.
\end{lemma}

Using Lemma~\ref{lemma1}, the PE condition may be checked online by determining whether or not $\xi$ has a minimal eigenvalue. During implementation, to maintain this condition, we introduce adequate exploration noise to the control as necessary.

Define next the estimation error for the value function as $\Tilde{V}_+ = V^*_+ - \hat{V}_+$. If the system state
$x_+(\tau)$ is bounded for a stabilizing control and $u_+(x)$, $x_+(\tau)$ and $\Delta \phi(x_+)$ are persistently excited, then we have the following theorem adopted from \cite{chen2019adaptive}.
\begin{theorem} \label{theorem1}
Consider system~\eqref{system} with the updating policy law~\eqref{weightupdate}, then: 
\begin{enumerate}
    \item If there is no neural network approximation error, i.e.
$\epsilon(x_+) = 0$, the error $\Tilde{V}_+$ will converge to zero in finite time $\tau_1 > 0$
\item In the presence of a neural network approximation error, $\Tilde{V}_+$ will converge to a small bounded set around its optimal control solution $V^*_+$ in finite time $\tau_1 > 0$
\end{enumerate}
\end{theorem}

Following the learning procedure described in this section, the closed-loop performance of the learned-control system, comprised of \eqref{system} and $u_\text{learned}=u_++u_-,$ will approximate that of the original control system, provided that $\varepsilon$ is sufficiently small or $T$ is sufficiently long. This will be further verified next through simulation examples.

%% file: Sections/5.Simulation.tex
\section{Examples} \label{simulation}
In this section, we demonstrate the effectiveness of our approach over three examples. 
\begin{figure*}
\centering
\begin{tabular}{ccc}
\includegraphics[width = 5cm]{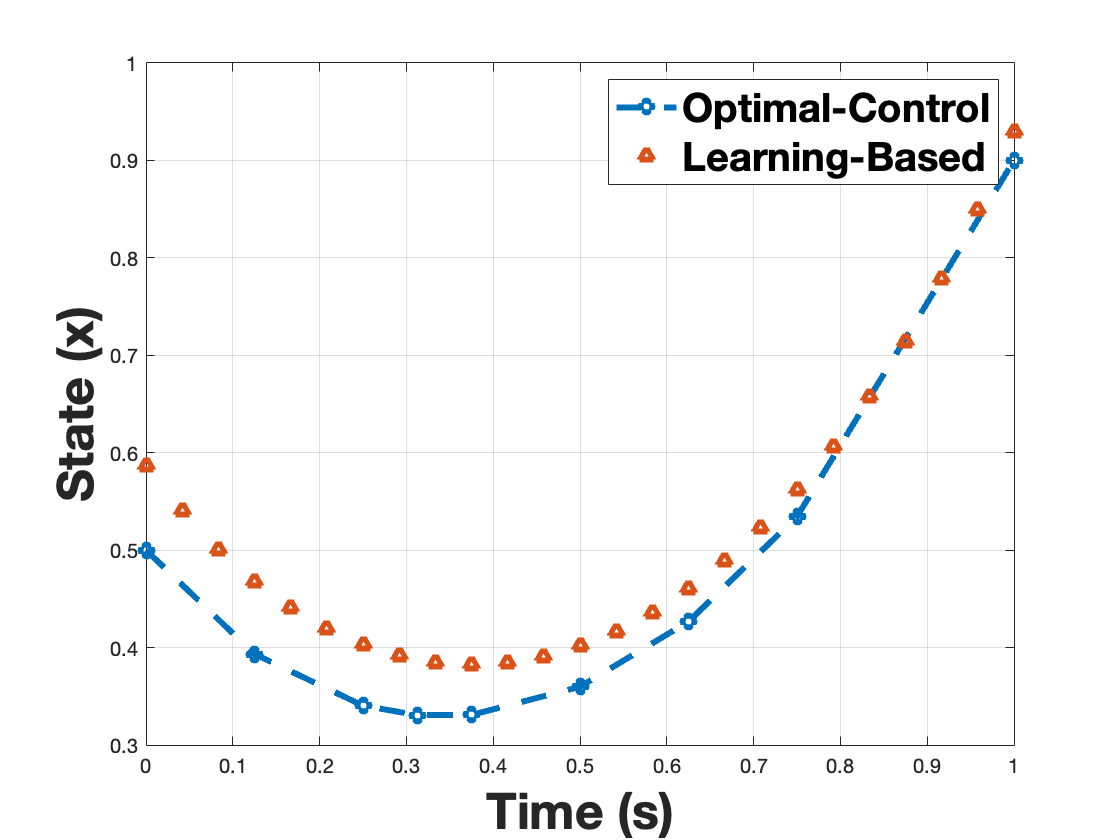} &
\includegraphics[width = 5cm]{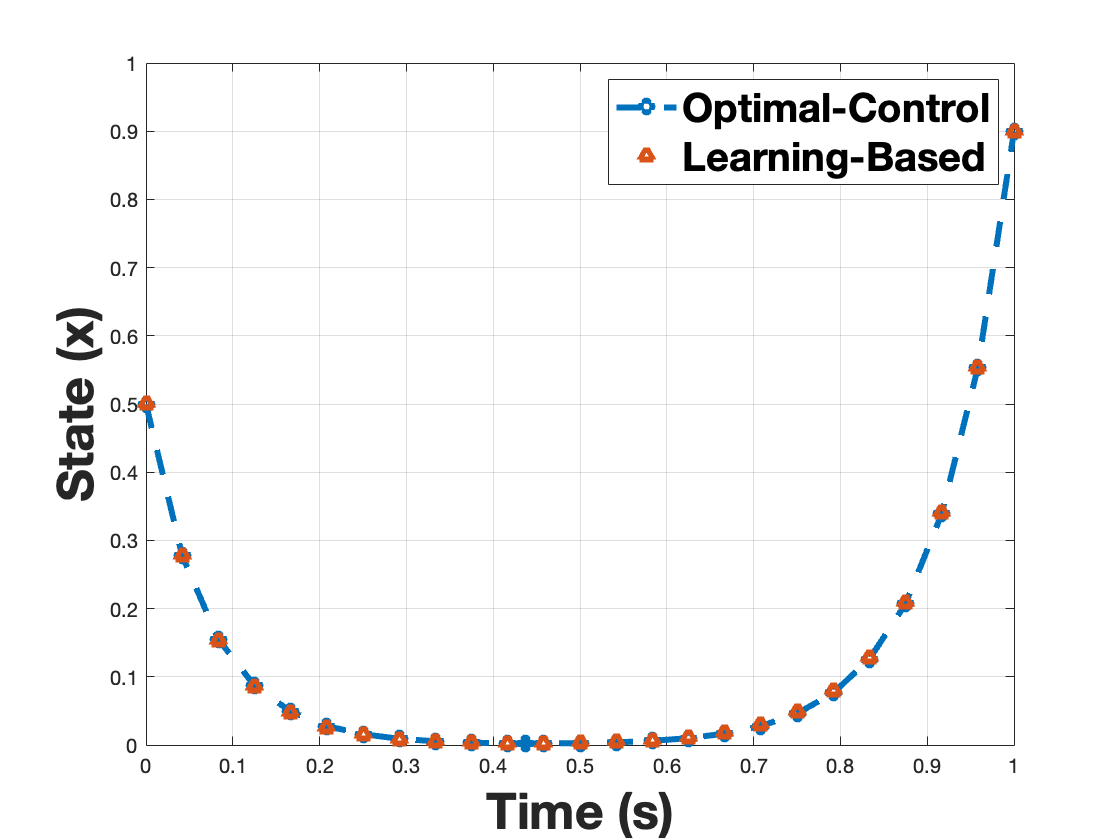} &
\includegraphics[width = 5cm]{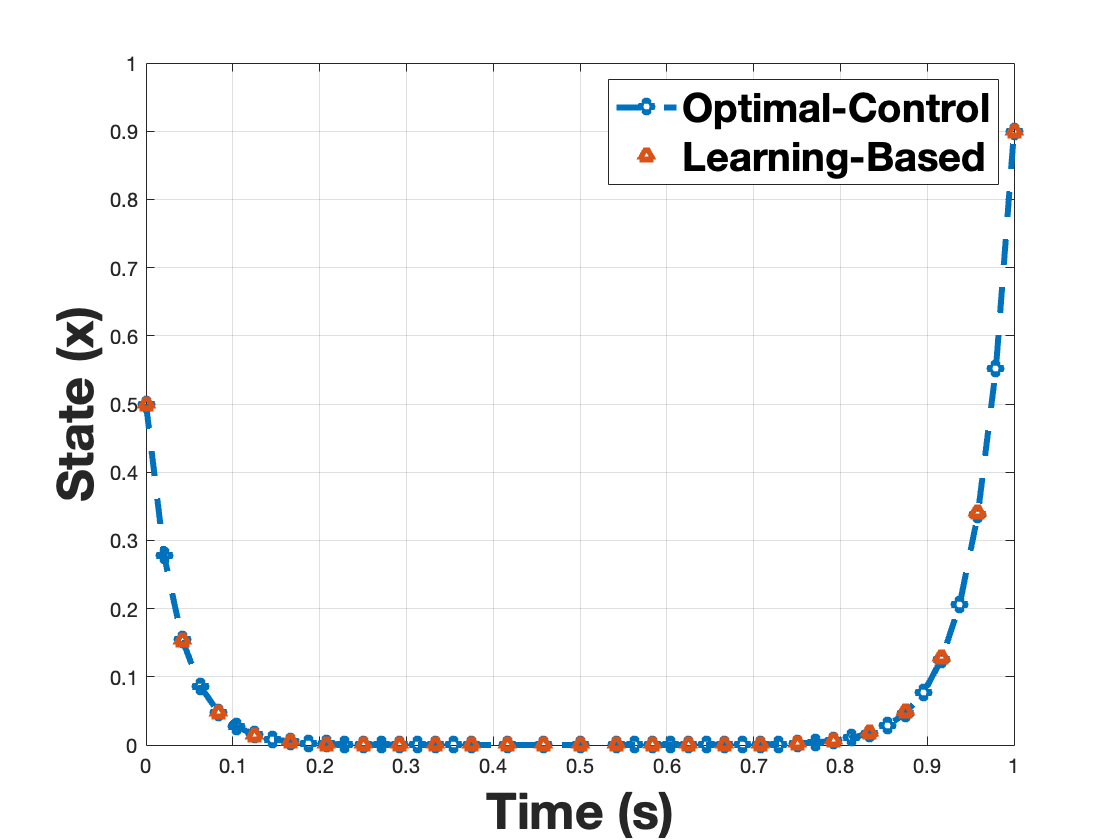} \\
\includegraphics[width = 5cm]{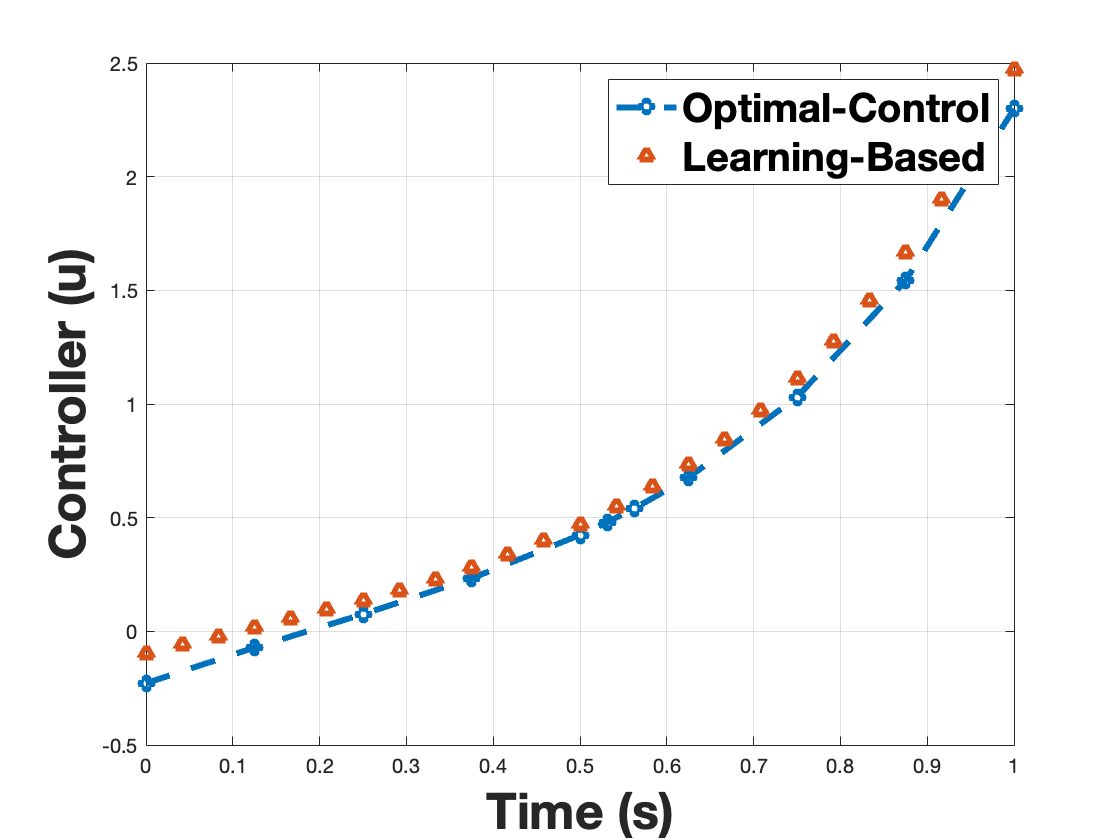} &
\includegraphics[width = 5cm]{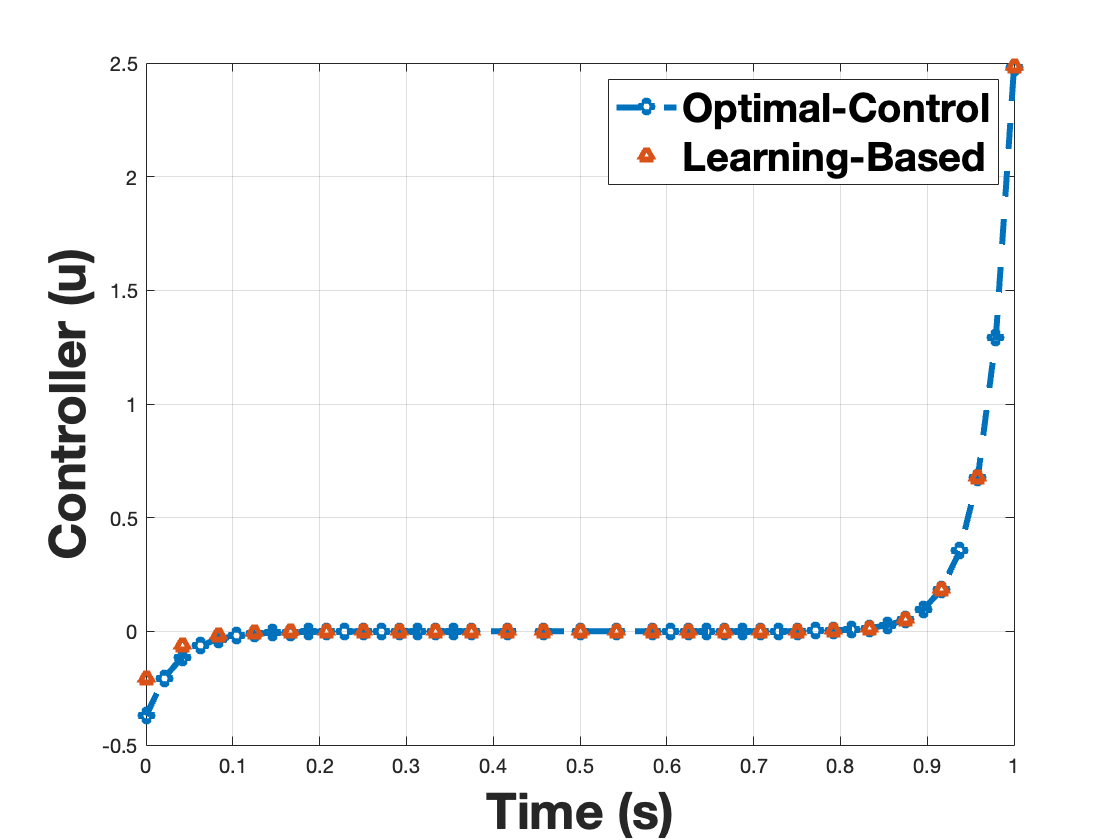} &
\includegraphics[width = 5cm]{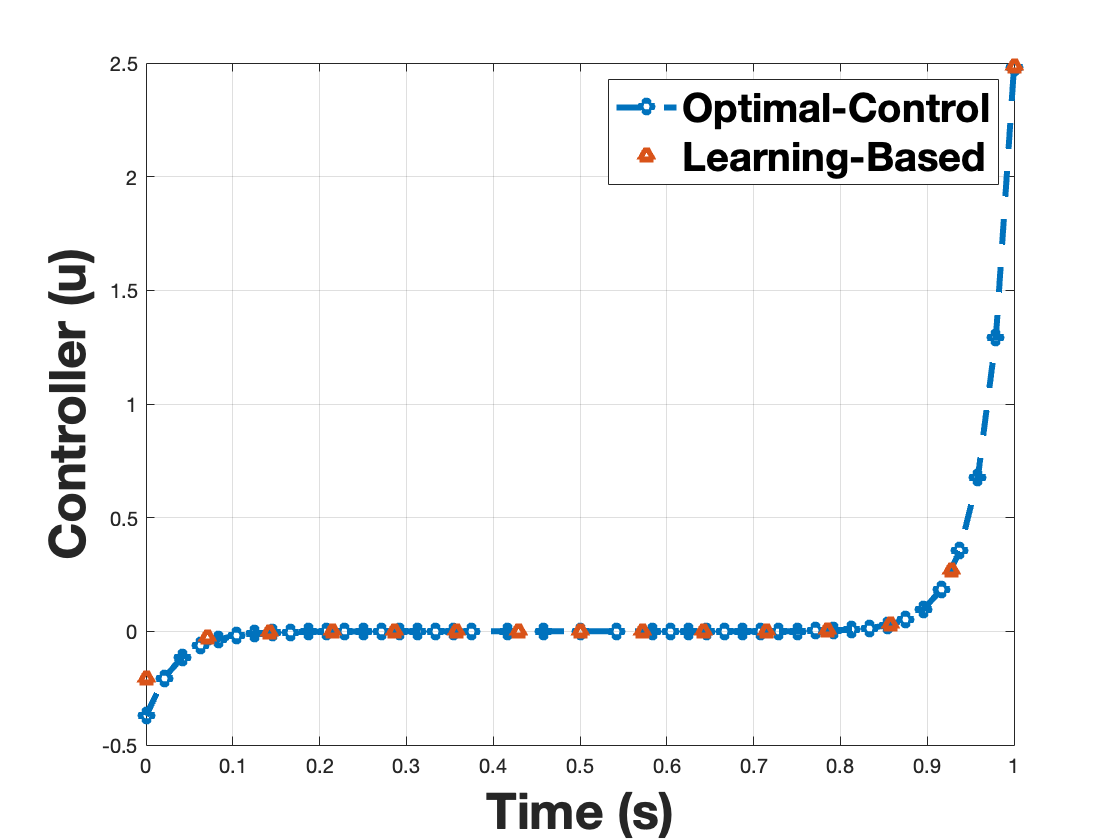} \\
\end{tabular}
\caption{(Top) The state space trajectory $x$ and (Bottom) the control law $u$ for the system~\eqref{objectdynami} is depicted using plots for various values of $\varepsilon = 0.5, 0.1$, and $0.05$, respectively.}
\label{1}
\end{figure*}
\begin{figure*}
\centering
\begin{tabular}{ccc}
\includegraphics[width = 5cm]{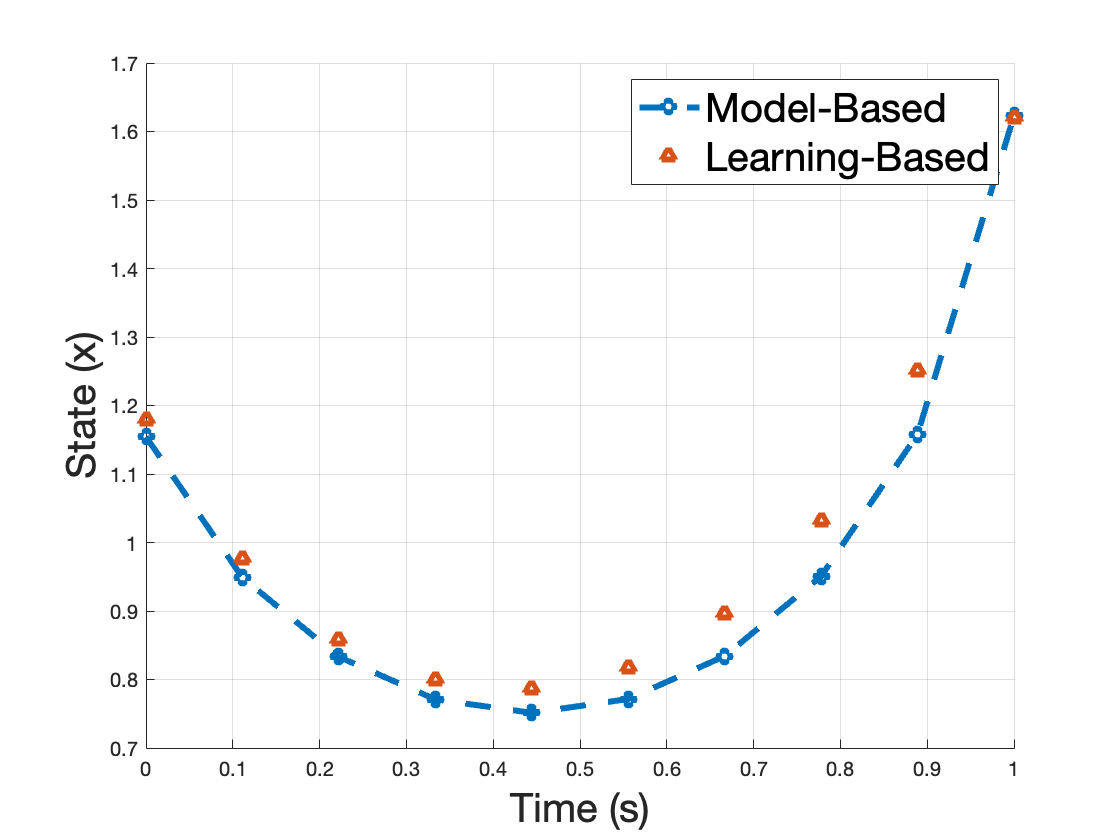}&
\includegraphics[width = 5cm]{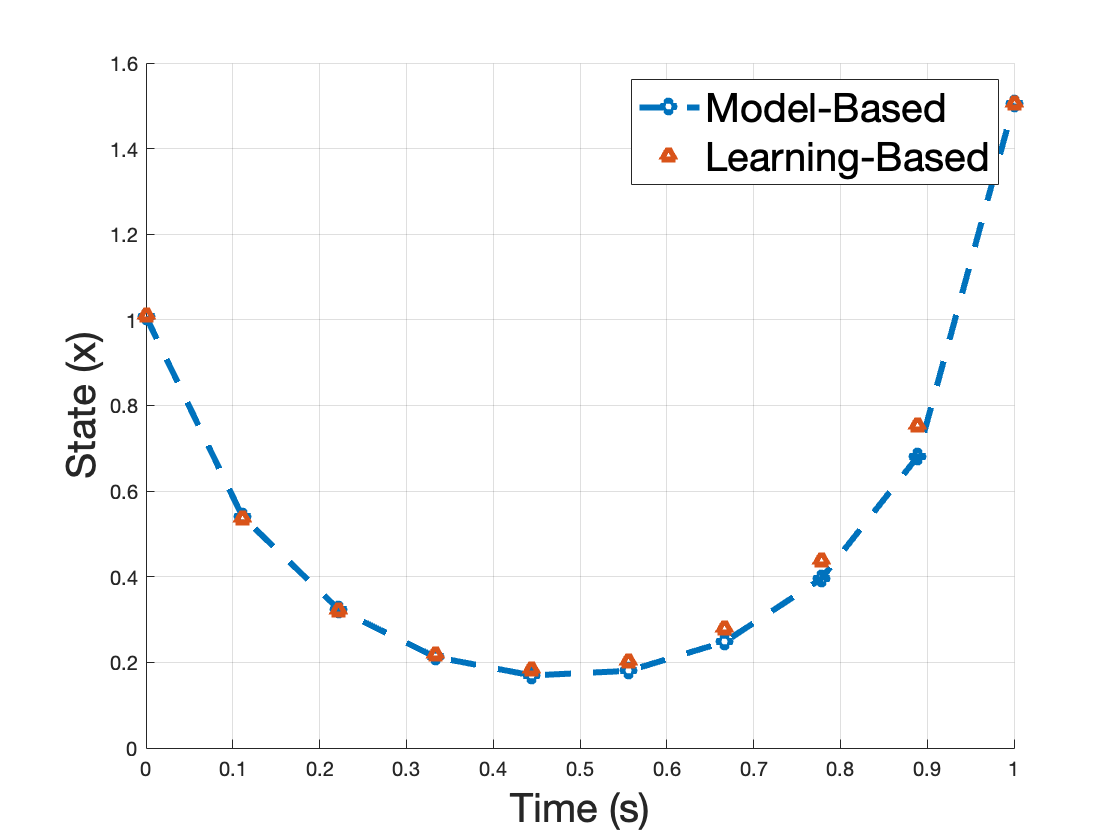}&
\includegraphics[width = 5cm]{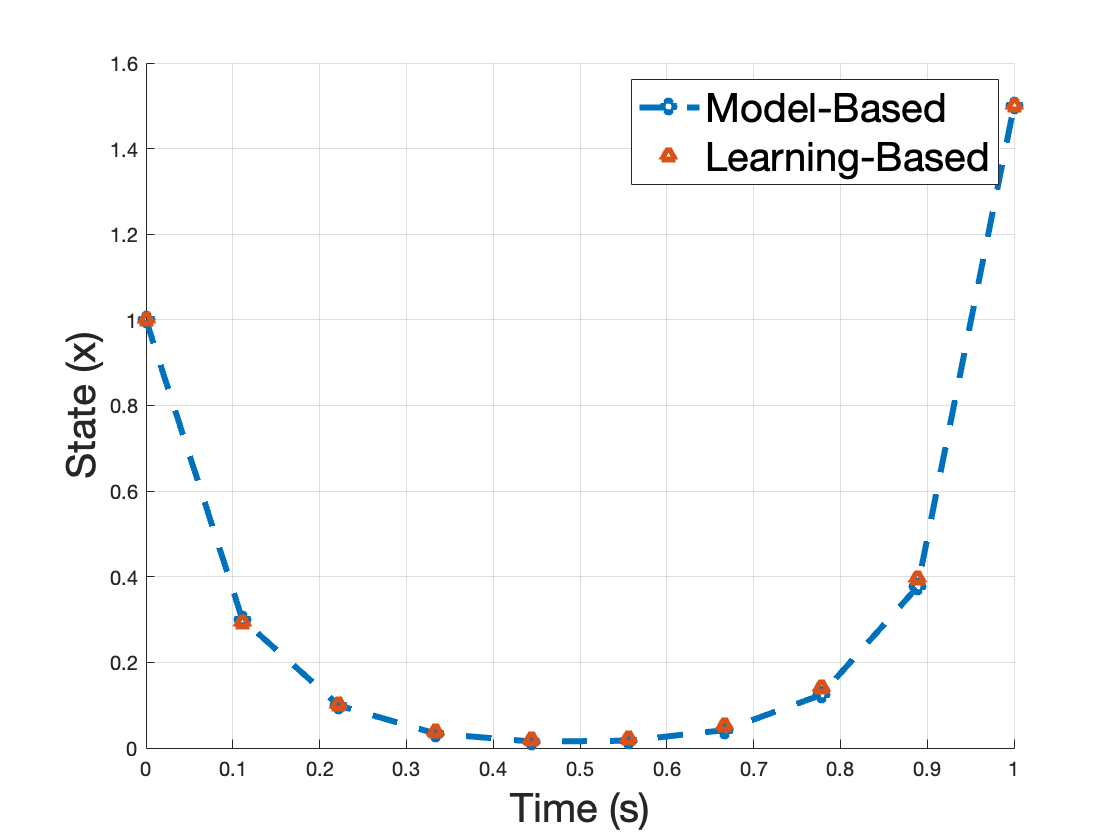} \\
\end{tabular}
\caption{The state space trajectory for the state $x$ for the system~\eqref{nonlinear} is depicted using plots for various values of $\varepsilon = 0.5, 0.2$, and $0.1$ accordingly.}
\label{3}
\end{figure*}
\subsection{RL Circuit}
Consider a circuit with the elements of resistance ($r$) and an inductor ($l$). The equation for the state space of the circuit is as follows:
\begin{equation} \label{objectdynami}
    \dot{x} = -\frac{r}{l}x + \frac{1}{L}u,
\end{equation}
 where $x$ is the state (the current through the circuit), and $u$ is the control input (circuit voltage). In this particular scenario, $f(x) = -\frac{r}{l}$ and $g(x) = \frac{1}{l}$. We will proceed under the assumption that the values of the inductor $L$ and resistor $R$ are not known. A controller $u$ is needed so that the current $x$ may be set to desired values at the beginning and end of the process while maintaining the minimum possible values for both the current and the voltage in between. To put it another way, the controller is necessary in order to achieve the optimal value of the objective function~\eqref{objectivefunction} while still ensuring that the system's initial and final conditions are $x_0 = 0.5$ and $x_T = 0.9$ respectively. The cost function parameters are $\mathcal{S}(x)=1$, $R=1$. For the sake of simulation, we will suppose that the inductor is described by the equation $l = 1 H,$ and $r = 1 \Omega$. In order to guarantee that the weights will eventually converge, we stimulate the system by means of the exploration noise signal $e = 2\sin t$. We choose the activation function vector for both forward and backward regulators as $\phi(x) = \frac{1}{2}x^2$, and the weights are initialized randomly between 0 and 1. After the training, following the algorithm~\ref{algo}, the weights converges as, $W_{+} = 0.41$, $W_{-} = -2.41$. \par 
Using the information provided by the functions $f(x)$ and $g(x)$, the two-value boundary problem solver available through the MATLAB command \textit{bvp4c} is used to identify the optimal state trajectory as well as the ideal controller. Next, we evaluate how well the state space trajectory and the controller our learning approach produced measure up against the characteristics of the best possible controller. \par
\begin{figure*}
\centering
\begin{tabular}{ccc}
\includegraphics[width = 5cm]{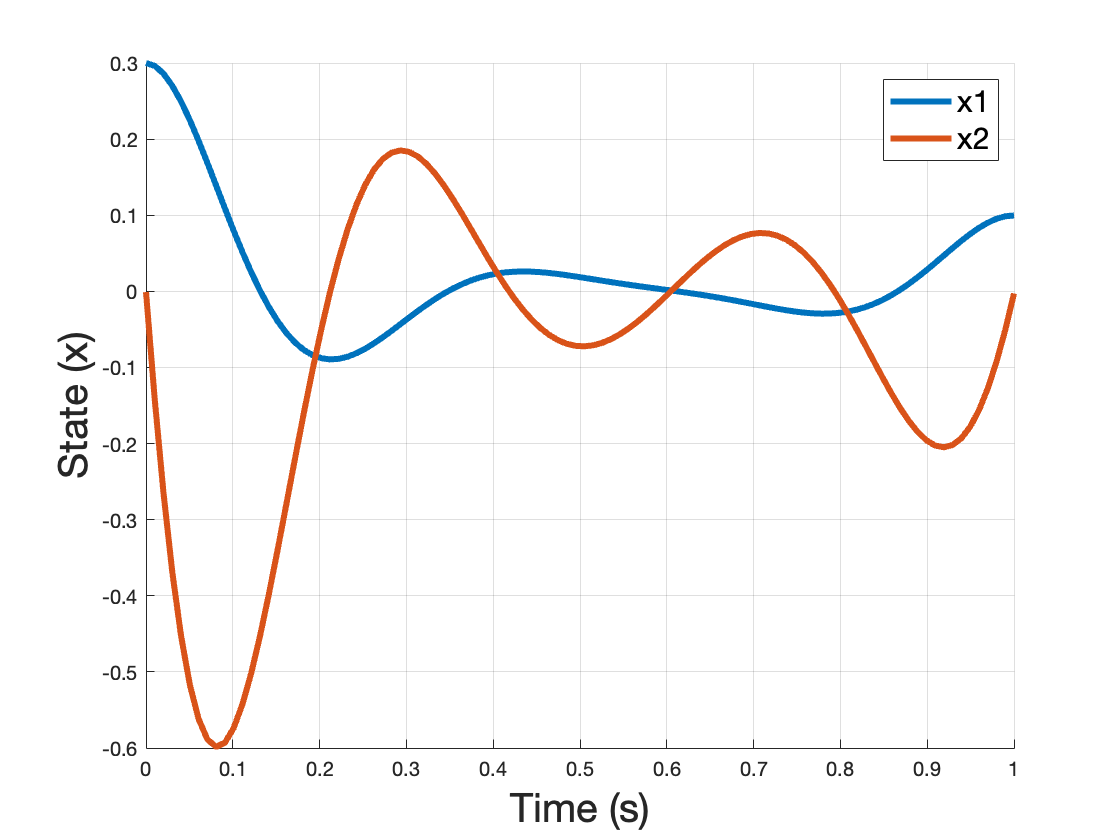} &
\includegraphics[width = 5cm]{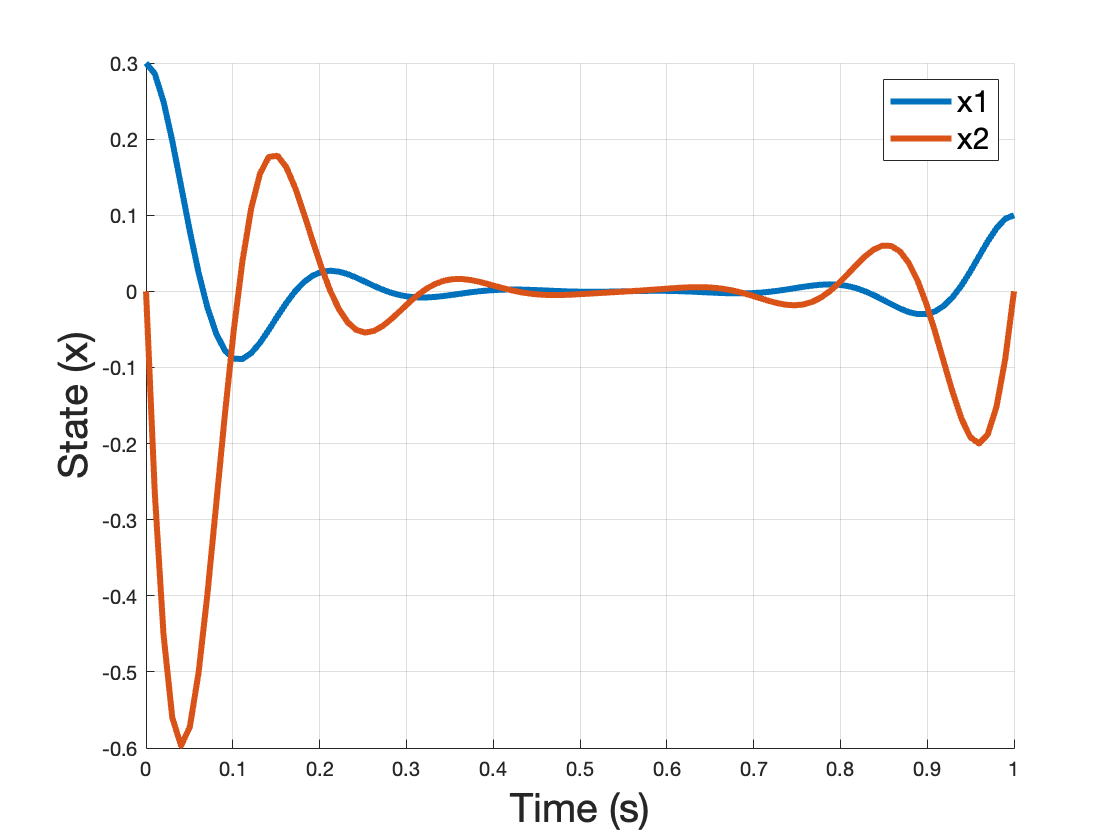}&
\includegraphics[width = 5cm]{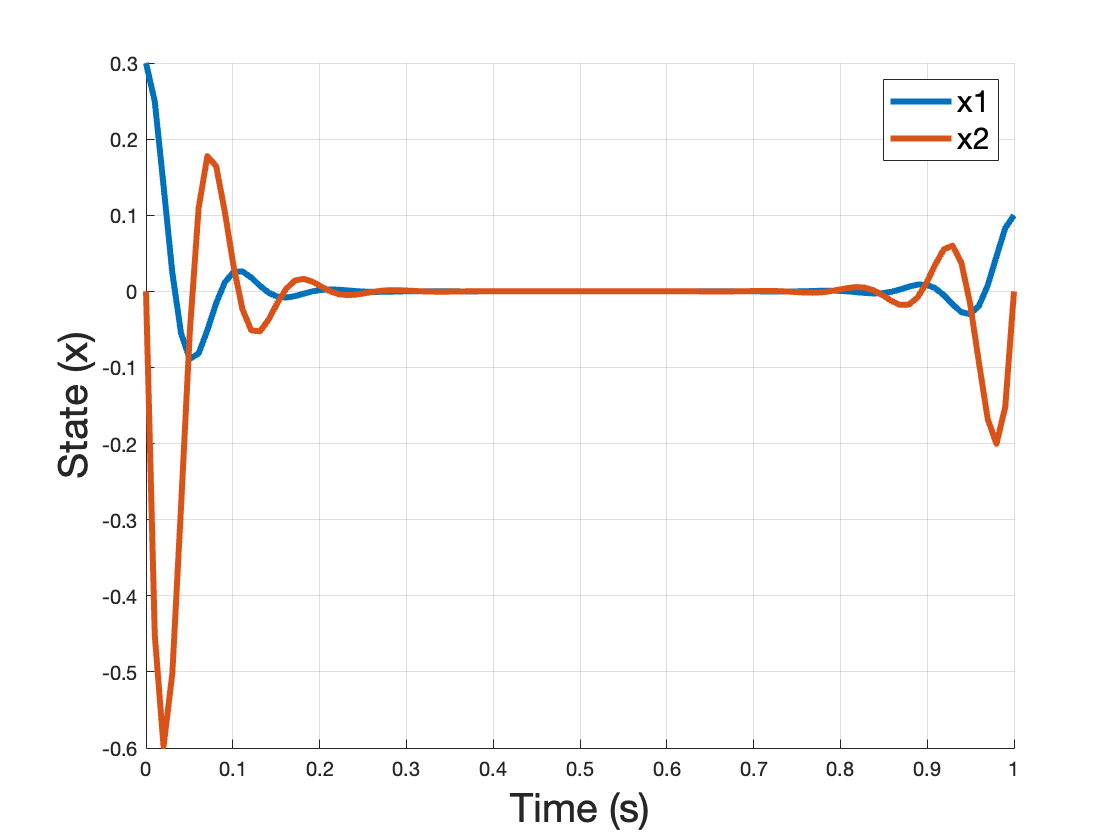} \\
\end{tabular}
\caption{The state space trajectory of the robotic manipulator system~\eqref{robot} is depicted using plots for various values of $T = 5,$ sec, $10$ sec, and $20$ sec, from left to right.}
\label{4}
\end{figure*}

As shown in Fig.~\ref{3} and Fig.~\ref{4}, as the value of $\varepsilon$ is decreased (or the control time period $T$ is increased), the learning-based controller becomes closer and closer to the optimal controller. This suggests that the learning controller has a less-than-ideal performance but will eventually converge to the optimal performance as the time interval $T$ becomes larger.
\subsection{Nonlinear System}
Consider a scalar nonlinear system
\begin{equation} \label{nonlinear}
    \dot{x} = x^3 + u 
\end{equation}
Our objective is to find a controller $u$ that steers the state $x$ from $x_0 =1$ to $x_T = 1.5$ while minimizing the objective function~\ref{objectivefunction} with the parameters, $\mathcal{S}(x) = 1$ and $R = 1$. In order to guarantee that the weights will eventually converge, we stimulate the system by means of the exploration noise signal $e = 2\sin t$. We choose the activation function vector for both forward and backward regulators as $\phi(x) = \left[x^2, x^4 \right]^\top$, and the weights are initialized randomly between -1 and 1. Following Algorithm~\ref{algo}, the weights converge as $W_{+} = \left[0.9740, 0.6933\right]^\top $, $W_{-} = \left[-0.9385, -0.6360\right]^\top$. \par 

In this example, our aim is to compare the proposed data-driven approach to that achieved by following the singular perturbation method described in Section \ref{method} had the system model been known. In this case, one can solve for $\frac{dV}{dx}$ by solving~\eqref{HJB+} and~\eqref{HJB-}~\cite{anderson1987optimal}. This results in the forward and backward value functions as: 
\begin{align} 
    V_{+}(x)&=\frac{1}{2} x^2\left(\sqrt{1+x^4}+x^2\right)+\frac{1}{2} \ln \left(\sqrt{1+x^2}+x^2\right),  \label{valueana+}\\
    V_{+}(x)&=-\frac{1}{2} x^2\left(\sqrt{1+x^4}+x^2\right)+\frac{1}{2} \ln \left(\sqrt{1+x^2}+x^2\right). \label{valueana-}
\end{align}
Substituting~\eqref{valueana+} in~\eqref{control+} and~\eqref{valueana-} in~\eqref{control-} yields the forward and backward regulators, respectively.

The learning-based approach and model-based controllers simulations are given in Fig.~\ref{3}. It can be seen that both controllers give similar results. 
\subsection{Robot Manipulator}
Consider a robotic manipulator, which is modeled by the equations~\cite{yu2020global}:
\begin{equation} \label{robot}
    \dot{x} = \left[\begin{array}{c}
x_2 \\
-2x_2-10\sin(x_1)
\end{array}\right] + \left[\begin{array}{c}
0 \\
1
\end{array}\right] u,
\end{equation}
where $x_1$ and $x_2$ are the angular and angular velocities, respectively and $u$ is the input torque. Let, $x_0 = [0.3 \; 0]^\top$ and $x_T = [0.1 \; 0]^\top$. The cost function parameters are specified as $Q=\left[\begin{array}{cc}
10 & 1 \\
1 & 10
\end{array}\right], \quad R=1$. A torque $u$ is needed to set the state $x$ at the desired boundary values while keeping the cost \eqref{objectivefunction} minimum. 

We follow Algorithm \ref{algo} to solve the problem and excite the system using the exploration noise signal $e = 2\sin t $ during the training. We choose the activation function vector for both forward and backward regulators as $\phi(x) = \left[
x_1^2, x_1x_2, x_2^2, x_1^3, x_1^2x_2, x_1x_2^2, x_2^3\right]^{\top}$ and the weights are initialized randomly between 0 and 1. After the training, the weights converge to the values $W_{+} = \left[2.7713, 0.1235, 0.2622, 0.0829, 0.0161, -0.0037, -0.0033
\right]^{\top}$, 
$\begin{aligned}
    W_{-} =& [-2.5970, -0.0688, -0.2507, -0.0108,\\ &-0.0039, -0.0001, 0.0011]^{\top}
\end{aligned}$

\par 
From Fig.~\ref{4}, we see that, as the time period in which the system is to be controlled is large, the transient dynamics dominate at the initial and terminal state boundaries. Moreover, during the time period $(0, T)$, the system stays close to zero to optimize the control objective. 
 \begin{remark} 
     It is important to mention that the simulations of the learning-based controllers involved training separate controllers for both the Forward regulator and Backward regulator. These controllers were then combined by overlapping their trajectories onto each other.
 \end{remark}

%% file: Sections/6.Conclusion.tex
\section{Conclusion} \label{conclusion}
We proposed an optimal controller design using reinforcement learning for two-point boundary nonlinear systems over finite-horizon time periods. The proposed design leverages the fast time scale occurring at the boundary conditions to avoid the need to solve the time-varying HJB equation. Furthermore, we design a learning-based control strategy that does not need knowledge of the system model. We show that the accuracy of the controller performance improves as the problem time horizon increases. We presented simulation results to support our claims using three examples. In the future, we plan to investigate the robustness of the proposed approach to noisy data and uncertain control input function.  